\documentstyle[12pt]{article}
\newcommand{\const}{\mathop{\rm const}\limits}

\newcommand{\Var}{\mathop{\rm Var}\limits}

\newcommand{\mod}{\mathop{\rm mod}\limits}

\newcommand{\Law}{\mathop{\rm Law}\limits}

\newcommand{\Div}{\mathop{\rm Div}\limits}

\newcommand{\grad}{\mathop{\rm grad}\limits}

\begin{document}

\begin{center}

{\bf  A SIMPLE MONTE CARLO METHOD FOR }  \\

\vspace{3mm}

{\bf SOLVING OF NAVIER-STOKES EQUATIONS.} \par

\vspace{4mm}

 $ {\bf E.Ostrovsky^a, \ \ L.Sirota^b } $ \\

\vspace{4mm}

$ ^a $ Corresponding Author. Department of Mathematics and computer science, Bar-Ilan University, 84105, Ramat Gan, Israel.\\

\vspace{3mm}

E-mail: \  eugostrovsky@list.ru \\
eugeneiostrovsky@gmail.com\\

\vspace{3mm}

$ ^b $  Department of Mathematics and computer science. Bar-Ilan University,
84105, Ramat Gan, Israel.\\

\vspace{3mm}

 E-mail: sirota3@bezeqint.net \\

\vspace{4mm}
                    {\bf Abstract.}\\

 \end{center}

 \vspace{4mm}

 We offer a simple   Monte-Carlo method for solving  of the multidimensional initial value and non-homogeneous
problem for the Navier-Stokes Equations in whole space when the initial function and right hand side
 belong to the correspondent Sobolev-Lebesgue-Riesz space.

\vspace{4mm}

 \par

\vspace{4mm}

  \par

\vspace{4mm}

\section{Introduction. Notations. Statement of problem.}

\vspace{3mm}

 The mild solution $ u = u(x,t) $ of a Navier-Stokes equation in the whole space $  x \in R^d $ during
its lifetime $ t \in [0,T], \ 0 < T = \const \le \infty $ may be represented as a limit as $ n \to \infty, n = 0,1,2, \ldots $
in appropriate space-time norm of the following recursion:

 $$
 u_{n+1}(x,t) = u_0(x,t) + G[u_n, u_n](x,t), n=1,2,\ldots \eqno(1.0)
 $$
where $ u_0(x,t) $ is the well-known solution of {\it linear} heat equation  with correspondent initial value and right-hand side
and $ G[u,v] $ is bilinear unbounded pseudo-differential operator, including space-time convolution, Riesz's transform,
 see \cite{Calderon1}, \cite{Cannone1}, \cite{Fujita1}, \cite{Giga1}, \cite{Giga2},  \cite{Kato1}, \cite{Kato2}, \cite{Vishik1} etc.

 We offer a simple depending trial Monte-Carlo method (stochastic modelling),
 without method based on solving of non-linear system of algebraic equations offered in \cite{Temam1}, for  multiple parametric
integrals computation emerging in (1.0)  such that  by optimizing of the proportion between the amount of random variables
with different degrees of integrals the speed of convergence of a {\it random approximation}
$ u_{n,N} $  to $ u_n $ as $ \ N \to \infty $ is optimal: for many important space-times  norms $ ||\cdot||, $
for example mixed Bochner's type anisotropic Lebesgue-Riesz $ L_{r,p}(R^d,[0,T]) $ norms or Lebesgue-continuous norm

$$
\forall n = 1,2,\ldots  \ \Rightarrow {\bf E} ||  u_{n,N}  - u_n|| \le K(n) \ N^{-1/2}, \eqno(1.0a)
$$
analogously to the work \cite{Grigorjeva1} devoted to the linear integral equations. Here $  N $ denotes the
{\it general amount} of elapsed random variables. The proof of this proposal used CLT and LIL in the considered spaces,
see \cite{Dudley1}, \cite{Gine1},  \cite{Ostrovsky1}, \cite{Ostrovsky207}, \cite{Ostrovsky210}, \cite{Ostrovsky211},
\cite{Ostrovsky212}, \cite{Vaart1} etc.\par

\vspace{3mm}

{\bf Detail description of problem.} \par

\vspace{3mm}

We  consider in this article the  initial value problem  for the multivariate Navier-Stokes (NS) equations

$$
\partial{u}_t - 0.5\Delta u + (u \cdot \nabla)u = \nabla P, \ x \in R^d, \ d \ge 3,   \ t > 0; \eqno(1.1)
$$

$$
\Div (u)  =  0, \ x \in R^d, \ t > 0; \eqno(1.2)
$$

$$
u(x,0) = a(x), \ x \in R^d. \eqno(1.3)
$$
 Here as ordinary

$$
x = (x_1,x_2,\ldots,x_k,\ldots,x_d) \in R^d ; \  ||x||: = \sqrt{\sum_{j=1}^d x_j^2},
$$

$$
\partial g(x)= \grad g(x) = \{  \partial g/\partial x_i \},
$$

and
$$
u = u(t) = u(t,\cdot)  = u(x,t)  = \{ u_1(x,t), u_2(x,t), \ldots,  u_d(x,t) \}
$$
denotes the (vector) velocity of fluid in the point $ x $ at the time $  t,  \ P $   is
represents the pressure. \par
  Equally:

 $$
\partial{u_i}/\partial t  = 0.5 \sum_{j=1}^d \partial^2_{x_j} u_i - \sum_{j=1}^d u_j \partial_{x_j} u_i +
 \partial_{x_i} P,
 $$

$$
\sum_{j=1}^d \partial_{x_j} u_j = 0, \ u(x,0) = a(x),
$$

$$
\Div u = \Div \vec{u} = \Div \{ u_1, u_2, \ldots, u_d \} = \sum_{k=1}^d \frac{\partial u_k }{\partial x_k} = 0
$$
in the sense of distributional derivatives.\par
 As long as

 $$
 P = \sum \sum_{j,k = 1}^d R_j R_k (u_j \cdot u_k),
 $$
where $ R_k = R_k^{(d)} $ is the $ k^{th} \ d \- $ dimensional Riesz transform:

$$
 \ R_k^{(d)}[f](x) = c(d) \lim_{\epsilon \to 0+} \int_{ ||y|| > \epsilon} ||y||^{-d} \Omega_k(y) \ f(x-y) \ dy,
$$

$$
c(d) = -\frac{\pi^{(d+1)/2}}{\Gamma \left( \frac{d+1}{2} \right) }, \ \Omega_k(x)= x_k /||x||,
$$
the system  (1.1) - (1.3) may be rewritten as follows:

$$
\partial{u}_t = 0.5 \Delta u - (u \cdot \nabla)u  +   Q \cdot \nabla \cdot (u \otimes u), \ x \in R^d, \ t > 0; \eqno(1.4)
$$

$$
\Div (u) = 0, \ x \in R^d, \ t > 0; \eqno(1.5)
$$

$$
u(x,0) = a(x), \ x \in R^d, \eqno(1.6)
$$
where   $ Q $ is multidimensional  Helmholtz-Weyl projection
operator, i.e., the $ d \times d  $ matrix pseudo-differential operator in $ R^d $ with the matrix symbol

$$
a_{i,j}(\xi) = \delta_{i,j} - \xi_i \xi_j /||\xi||^2, \hspace{5mm} \delta_{i,j} = 1, i = j; \delta_{i,j} = 0, \ i \ne j.
$$

\vspace{4mm}
 The consistent regularization of the Riesz transform looks as follows:

$$
  R_{k, \epsilon}^{(d)}[f](x) = c(d)  \int_{R^d } [ \epsilon + ||y||^{-d}] \ \Omega_k(y) \ f(x-y) \ dy,
$$
herewith

$$
|| R_{k, \epsilon}^{(d)}[f] - R_{k}^{(d)}[f]||_p \le C(p) \ \epsilon \ ||f||_p,
$$
see  \cite{Stein1}, chapters 4,5. \par

 Note that the last representation of the Riesz's potential may be used by its Monte Carlo  computation, if we will use
the polar coordinates and the density of applied random variables to be proportional to the kernel $  \epsilon + ||y||^{-d}. $ \par

\vspace{4mm}
 At the same considerations may be provided for the NS equations with external density of  force $  f = f(x,t): $

$$
\partial{u}_t = 0.5\Delta u - (u \cdot \nabla)u  +   Q \cdot \nabla \cdot (u \otimes u) + f(x,t), \ x \in R^d, \ t > 0;
$$

$$
u(x,0) = a(x), \ x \in R^d.
$$
see  \cite{Giga1} -  \cite{Giga4}, \cite{Koch1},  \cite{Masuda1}, \cite{Ogawa1}, \cite{Temam1}.\par
 More detail, the considered here problem may be rewritten as follows:

$$
u(x,t) = e^{0.5 t \Delta} a(x)  + G [u,u](t) \stackrel{def}{=} u_0(x,t) + G [u,u](t) + v[f](x,t),
$$
where

$$
v[f](x,t) = v(x,t) = v = \int_0^t ds \int_{R^d} w_{t-s}(x-y) \ f(y,s) \ dy =
$$

$$
\int_0^t  w_{t-s}(\cdot)*f(\cdot,s) \ ds.
$$

 We will denote further for simplicity

 $$
 \int f(y) dy = \int f  = \int_{R^d} f(y) dy,
 $$

 $$
 w \odot u (x,t) \stackrel{def}{=} \int_0^t \int_{R^d} w_{t-s}(x-y) u(y,s) ds dy  =
  \int_0^t \int w_{t-s}(x-y) u(y,s) ds dy.
 $$

\vspace{4mm}
 {\it We will understand henceforth as a capacity of the solution (1.4)-(1.6) the vector-function  } $ u = \vec{u} =
 \{ u_1(x,t),  u_2(x,t), \ldots,  u_d(x,t) \} $ {\it  the so-called mild solution,} see \cite{Miura1}.  \par
\vspace{4mm}

 Namely, the vector- function  $ u = u(\cdot,t) = u(x,t)  $ satisfies almost everywhere in the time $ t $ the following {\it non-linear
 integral-differential equation: }

$$
u(t) = e^{0.5 t \Delta} a + \int_0^t e^{0.5(t-s)\Delta } [ (u \cdot \nabla)u(s)  +   Q \cdot \nabla \cdot (u \otimes u)(s) ] ds \stackrel{def}{=}
$$

$$
e^{0.5 t \Delta} a  + G [u,u](t) \stackrel{def}{=} u_0(x,t) + G [u,u](t), \eqno(1.7)
$$
 the operator $  \exp( 0.5 t \Delta) $ is the classical convolution  integral operator with heat kernel: \par

$$
u_0(x,t) :=  e^{0.5 t \Delta}[ a](x,t) = w_t(x)*a(x),
$$
where $ G(u,u)  \stackrel{def}{=}  F(u,u) = F(u), $

$$
F(u,v) = \int_0^t e^{0.5 (t-s)\Delta } B[u,v](x,s) \ ds =  \int_0^t \ \int_{R^d} \ w_{t-s}(x-y) \ B[u,v](y,s) ds \ dy, \eqno(1.8)
$$

$$
 B(u,v):=  (u \cdot \nabla)v(s)  +   Q \cdot \nabla \cdot (u \otimes v)(s),
$$

$$
w_t(x) \stackrel{def}{=} (2 \pi t)^{-d/2}  \exp \left( - \frac{|| x ||^2}{2 t} \right).
$$
The convolution between two functions   $ r = f(t), \ g(t)  $ defined on the set $ R_+ = (0,\infty)  $ is defined as usually

$$
f \odot g(t) := \int_0^t f(t-s) \ g(s) \ ds  = g \odot f(t)
$$
("time-wise" convolution) and between two, of course, measurable  vector-functions  $ u(x), v(x) $ defined on the whole space $ x \in R^d $

$$
u*v(x) = \int_{R^d} u(x-y) \ v(y) \ dy,
$$
("space-wise" or  "coordinate-wise" convolution).\par
  The authors hope that this notations does not follow the confusion. \par

\vspace{3mm}

   More results about the existence, uniqueness, numerical methods, and a priory estimates in the different Banach function spaces:
 Lebesgue-Riesz $  L_p, $ Morrey, Besov for this solutions see, e.g. in  \cite{Cui1}- \cite{Zhang1}. The first and besides famous
 result belong to J.Leray  \cite{Leray1}; it is established there in particular the {\it global in time}
 solvability and  uniqueness of NS system in the space $ L_2(R^d)  $ and was obtained  a very interest a priory estimate for solution.\par
  The quantitative estimations for solution and lifespan of NS equations  in some rearrangement invariant spaces see, e.g. in
 \cite{Ostrovsky201}, \cite{Ostrovsky202}. \par

   T.Kato in \cite{Kato1} proved in particular that if the initial function $ a = a(x)  $ belongs to the space $  L_d(R^d) $
(in our notations), then there exists a positive time value $ T > 0  $ (lifespan of solution) such that the solution of NS system
$ u = u(x,t) $ there exists for $ t \in (0,T), $  is smooth and satisfy some a priory integral estimates.\par
 Furthermore,  if the norm $ ||a||L_d(R^d)  $ is sufficiently small, then $  T = \infty, $ i.e. the solution $ u = u(x,t) $
is global.\par
 The {\it upper} estimate for the value $ T, $ conditions for finite-time blow-up
 and asymptotical behavior of solution as $ t \to T - 0 $ see in the articles \cite{Ball1}, \cite{Benameur1},
\cite{Chae1},  \cite{Cui1}, \cite{Gallagher1}, \cite{Montgomery1},  \cite{Seregin1}, \cite{Seregin2} etc.\par

\vspace{4mm}

 With regards the numerical methods for solving of NS equations, we note only the classical monograph \cite{Temam1} and
 articles  \cite{German1}, \cite{Vishik1}, \cite{Serrin1}. \par

\vspace{4mm}

{\bf  Our purpose in this short report is to offer some modification of the optimal Monte - Carlo  method for solving of NS equations during
 the lifespan of solution } $  T. $ \par

\vspace{4mm}

 The essence of the proposed method can be explained very simply: we will write the approximation $  u_n(x,t) $ in the {\it  explicit  view}
through multiple sums of multiple  parametric integrals of increasing dimension, to calculate which may be used the so-called depending
trial Monte Carlo method. \par

\vspace{3mm}

 With regard to the Monte-Carlo method for the {\it function computation } (in the other terms, depending trials method) see
\cite{Frolov1}, \cite{Grigorjeva1}, \cite{Ostrovsky204},  \cite{Ostrovsky207}.\par

 Note that Monte-Carlo  method can not give a very high precision, but it is very simple. For instance, it does not use the
solving of system of (non-linear!)  algebraic equations and following is stable (robust). \par

 \vspace{3mm}

\section{ Some Notations, with Clarification. The essence of the method.} \par

\vspace{3mm}

 As ordinary, for the measurable function $ x \to u(x), \ x \in R^d $

 $$
||u||_p = \left[ \int_{R^d} |u(x)|^p \ dx  \right]^{1/p}. \eqno(2.1)
 $$

{\it  Multidimensional case.} \par
 Let $ u = \vec {u} = \{ u_1(x), u_2(x), \ldots, u_d(x)  \}  $ be measurable vector-function: $ u_k: R^d \to R. $ We can define
as ordinary the $ L_p, \ p \ge 1 $  norm of the function $ u $ by the following way:

$$
||u||_p := \max_{k=1,2, \ldots,d} || u_k||_p, \ p \ge 1. \eqno(2.2)
$$

\vspace{4mm}

We present now some important results belonging to   T.Kato \cite{Kato1}; see also the article of
 H.Fujita  and T.Kato  \cite{Fujita1}.  Let us consider the following recursion:

$$
u_{n+1}(x,t) = u_0(x,t)  + G [u_n, u_n](t)= u_0(x,t) + G [u_n, u_n](t) \eqno(2.3)
$$

with initial condition for iterations

$$
u_0(x,t) = (2 \pi t)^{-d/2} \int_{R^d} \exp \left( - \frac{||x-y||^2}{2t}  \right) \ a(y) \ dy +
$$

$$
\int_0^t ds \int_{R^d} (2 \pi t)^{-d/2} \exp \left( - \frac{ ||x-y||^2}{2(t-s)}  \right) \ f(y,s) \ dy.
$$

 So, the recursion (2.3) may be rewritten as follows:

$$
u_{n+1}(x,t)= u_0(x,t) +
(2 \pi t)^{-d/2} \int_0^t \ ds \ \int_{R^d} \exp \left( - \frac{||x-y||^2}{2(t-s)}  \right)  \ F[u_n](y,s) \ dy. \eqno(2.4)
$$

 H.Fujita  and T.Kato  \cite{Fujita1}, \cite{Kato1}  proved that if $ a(\cdot) \in L_d(R^d), \ \Div a = 0, $
 then there exists a positive value $ T = T(||a||_d)  $  (lifespan of solution) such that the iteration sequence
$  u_n(\cdot, \cdot) $ converges in the senses (1.9) to the uniquely defined solution of NS equations.\par

 Furthermore, if the norm $ ||a||_d $  is sufficiently small, then $ T = \infty $ (global solution). \par
 The quantitative lower bound for lifespan $  T  $ and some quantitative a priory estimation in Grand Lebesgue Spaces (GLS)
for solution $ u = u(x,t) $ see, e.g. in \cite{Ostrovsky201}, \cite{Ostrovsky202}.\par

Moreover, for all the values $ \delta \in (0,1) $
there exist a constants $ q_1 = q_1(\delta), q_2 \in (0,1) $  such that

$$
|| u_n - u||B( (0,T), L_{d/\delta} )  \le C_1(a; d,\delta) \ q_1^n
$$
and

$$
|| \nabla u_n - \nabla u||B( (0,T), L_{d} )  \le C_2(a; d) \ q_2^n.
$$

\vspace{4mm}

{\bf A.} Note that if $  f(x,t) = 0 $ then

 $$
 u_0(x,t) = u_0(x,t) = (2 \pi)^{-d/2} \int_{R^d} \exp \left( - \frac{||z||^2}{2} \  \right) \ a(x + z \sqrt{t}) \ dz =
 $$

$$
{\bf E} a(x + \xi \cdot \sqrt{ t}),  \eqno(2.5)
$$
where the {\it random vector}  $  \xi $ has a  {\it Gaussian } $ d \ - $ dimensional standard  distribution. \par
 This imply that the random vector  $  \xi $ has a  {\it Gaussian } $ d \ - $ dimensional  distribution with parameters

 $$
 {\bf E} \xi = 0,  \hspace{5mm} \Var \xi = I_d
 $$
be an unit matrix of a size $  d \times d. $ \par

  Let $ \{  \xi_i \}, \ i = 1,2,\ldots, N  $  be a sequence of independent Gaussian  $ d \ - $ dimensional standard  distributed
 random vectors.  The Monte-Carlo approximation  $ u_{0,N} =u_{0,N}(x,t) $ in that its modification which is called "depending trial method"
 \cite{Frolov1} for $ u_0(x,t) $ has a view

 $$
u_{0,N}(x,t)= u_{0,N}[a](x,t):= N^{-1} \sum_{i=1}^N a(x + \xi_i \ \cdot \sqrt{t}).  \eqno(2.6)
 $$

\vspace{3mm}

{\bf B.} Analogously,  let us  consider the  $  d \ - $ dimensional heat equation with zero initial values but with external force:

$$
\partial_t u^{(0)}(x,t) = 0.5\Delta u^{(0)}(x,t)  + f(x,t), \  u^{(0)}(x,0) = 0. \eqno(2.7)
$$

 Then $ u^{(0)}(x,t) = u^{(0)}[f](x,t) =  $

 $$
  \int_0^t [w_{t-s} * f](x,s) \ ds  = \int_0^t ds \ \int_{R^d} (2 \pi (t-s))^{-d/2} \
 \exp \left(  - \frac{||x-y||^2}{2(t-s)} \right) \ f(y,s) \ dy =
 $$

 $$
 t \int_0^1 dv \ \int_{R^d} (2 \pi)^{-d/2} \ e^{ - ||z||^2/2} \ f(x + z \sqrt{t} \ \sqrt{(1-v) }, \ t \cdot v ) \ dz =
 $$

 $$
 {\bf E} \ [t \ f(x + \eta \sqrt{t} \ \sqrt{(1-\tau) }, \ t \cdot \tau)], \eqno(2.8)
 $$
 where the random {\it vector} $ \eta $  has a  Gaussian  $ d \ - $ dimensional standard  distribution,
 the random {\it variable} $ \tau $ is uniformly distributed on the unit segment $ [0,1]. $\par

 Let $ \{  \eta_i \}, \ i = 1,2,\ldots, N  $  be a sequence of independent Gaussian  $ d \ - $ dimensional standard  distributed
 random vectors and   $ \{\tau_i\},  \ i = 1,2,\ldots, N $ be a sequence of independent and independent on the $\{  \eta_i \}  $
 uniform distributed on the segment $ [0,1] $ random variables. \par

  The Monte-Carlo approximation  $ u^{(0)}_{N} =u^{(0)}_{N}(x,t) = u^{(0)}_{N}[f](x,t) $ in that its modification which is
called "depending trial method" for $ u^{(0)}(x,t) $ has a view

 $$
 u^{(0)}_{N}(x,t) = u^{(0)}_{N}[f](x,t) :=
\frac{t}{N} \ \left [ \sum_{i=1}^N f(x + \eta_i \sqrt{t} \ \sqrt{(1-\tau_i) }, \ t \cdot \tau_i) \right]. \eqno(2.9)
 $$

\vspace{3mm}

Of course, the "initial" function $ u_0(x,t) $ may be computed by means  of deterministic methods: finite differences, finite elements
method etc., as well as  the Riesz's transform computation. \par

\vspace{3mm}

{\bf C.} The expression for the  member $ G[u] $ or for the $ G[u_n] $ is  alike to the one in the formula (2.8) with replacing
$  f :=F = F[u](x,t). $\par

{\it Let us define the following iteration sequence:}

 $$
 u_{n+1, N(n+1)}(x,t) = u_{0,N(0)}[a](x,t) +
 $$

 $$
 \frac{t}{N(n+1)} \
 \left [ \sum_{i=1}^{N(n+1)} F[u_{ n, N(n)}]
 \left(x + \eta_i^{(n+1)} \sqrt{t} \ \sqrt{(1-\tau_i^{(n+1)}) }, \ t \cdot \tau_i^{(n+1)} \right) \right]. \eqno(2.10)
 $$

 Here   $ n = 0,1,2, \ldots,L; \  $ the number  $ L $ is the total number of  iterations,

$$
  (N(0), N(1), N(2), \ldots, N(L) ) =  (N^{(L)}(0), N^{(L)}(1), N^{(L)}(2), \ldots, N^{(L)}(L) )
$$
 is the sequence of integer numerical vectors ever-increasing dimension $ L+1, \ L = 1,2, \ldots; \ L \to \infty  $ such that

$$
\forall n = 0,1,2,\ldots L \ \Rightarrow \ \lim_{L \to \infty}  N^{(L)}(n)  = \infty.
$$

 The random vectors $ \eta_i^{(n+1)} $ have  $  d - $ dimensional standard Gaussian distribution, the random variables  $ \tau_i^{(n+1)} $
are uniformly distributed in the unit segment $ [0,1] $ and {\it all the introduced random vectors and  variables are independent. } \par

 The total number of spent random {\it variables}, i.e. including the vector coordinates $  M = M(L) $  may be calculated as follows:

 $$
 M = M(L) = 2N^{(L)}(0) + \sum_{d=1}^L (d+1) N^{(L)}(d). \eqno(2.11)
 $$

 But it is very hard to error estimate of this procedure, especially in important Banach functional norms for the computated function. \par

\vspace{4mm}

{\bf A. Let us consider now the  alternative approach.}  \par

\vspace{3mm}

Namely,  we denote

$$
 J_m(h)  = \partial w \odot \partial w \odot \ldots \odot \partial w \odot h =
$$

$$
\int_0^t ds_1 \int_0^{s_1} ds_2 \ldots \int_0^{s_{m-1}} ds_m \int_{R^{d m}} (2 \pi)^{-dm/2} (t-s_1)^{-d/2} (s_1 - s_2)^{-d/2} \ldots
(s_{m-1} - s_m)^{-d/2} \times
$$

$$
\exp \left[ - \frac{||x-y_1||^2}{2(t-s_1)}
- \frac{||y_1 - y_2||^2}{2(s_1 - s_2)} -  \ldots -  \frac{||y_{m-1} - y_m||^2}{2(s_{m-1} - s_m)}  \right] \times
$$

$$
 \left( \frac{x-y_1}{t-s_1} \cdot \frac{y_1 - y_2}{s_1 - s_2} \cdot  \ldots  \cdot \frac{y_{m-1} - y_m}{s_{m-1} - s_m}  \right) \ h(y,s) \ dy.\eqno(2.12)
$$

 We make the change of variables in the interior integral as follows:

 $$
 y_1 = x + z_1 \sqrt{t - s_1}, \ y_2 = y_1 + z_2 \sqrt{s_1 - s_2}, \ldots, y_m = y_{m-1} + z_m \sqrt{s_{m-1} - s_m},
 $$
with Jacobian

$$
\left[(t - s_1) (s_1 - s_2) \ldots (s_{m-1} - s_m)  \right]^{d/2},
$$
and after this - the substitution $  s_j = t \cdot \tau_j  $ with Jacobian $ t^{m/2}. $  The resulting transform may be written briefly as follows

$$
(y,s) = Y(z,\tau) = Y_{x,t}(z,\tau), \ y,z \in (R^d)^m = R^{dm}, \  s, \tau  \in (R_+)^m,
$$
where the values $ x \in R^d, \ t \in R^1_+  $ be a fix (temporarily).\par

We obtain:
 $  J_m(h) \cdot t^{-m/2} (2 \pi)^{dm/2} = $

$$
\int_{S(m)}  d \tau \int_{R^{d m}} \ dz  \frac{1}{(1-\tau_1)^{1/2}} \ \frac{1}{(\tau_1 - \tau_2)^{1/2}} \ldots
 \frac{1}{  (\tau_{m-1} - \tau_m)^{1/2} } \times
$$

$$
\exp \left[ - \frac{||z_1||^2}{2}
- \frac{||z_2||^2}{2} -  \ldots -  \frac{||z_m||^2}{2}  \right] \cdot \tilde{h}(z,\tau), \eqno(2.13)
$$
where  $  S(m)  $ is an  $ m \ - $  dimensional unit polygon  (simplex) of a form $  S(m) =  $

$$
 \{ \tau = \vec{\tau} = (\tau_1, \tau_2, \ldots, \tau_m):
    0 < \tau_1 < 1, 0 < \tau_2< \tau_1, 0 < \tau_3 < \tau_2, \ldots, 0 < \tau_m < \tau_{m-1}    \},
$$

and
$$
\tilde{h}(z,\tau) = \tilde{h}_{x,t}(z,\tau) = \prod_{j=1}^m z_j \cdot h( Y(z,\tau) ) =  \prod_{j=1}^m z_j \cdot h( Y_{x,t}(z,\tau) ).
$$

\vspace{3mm}

 It is easy to calculate

$$
\int \int \ldots \int_{S(m)}\frac{ds_1 ds_2 \ldots ds_m}{s_1^{\alpha_1} (s_2 - s_1)^{\alpha_2} (s_3 - s_2)^{\alpha_3} \ldots (s_m - s_{m-1})^{\alpha_m} } =
$$

$$
\frac{\prod_{k=1}^m \Gamma(1-\alpha_k)}{\Gamma(1 + \sum_{k=1}^m (1 - \alpha_k))},\ 0 \le \alpha_k < 1.
$$

 Denote  $ W_m =  W(m) = \pi^{m/2}/\Gamma(1 + m/2) \ - $ the volume of an unit ball of the Euclidean space  $  R^m, $

 $$
H_m(s) = \left[(1 - s_1)(s_1 - s_2) (s_2 - s_3) \ldots (s_{n-1} - s_n ) \right]^{-1/2}/W_m  \stackrel{def}{=}
$$

$$
R_m(s)/W_m, \ s \in S(m). \eqno(2.14)
 $$

 Evidently, $ \lim_{m \to \infty} W_m = 0.  $\par

 Recall in addition to the article \cite{Ostrovsky208} that the  generating function for the sequence

 $$
  W_m(\beta)  = \frac{1}{\Gamma(1 + n \beta)},
 $$
i.e. the function

$$
ML_{\beta}(z) = \sum_{n=0}^{\infty} \frac{z^n}{\Gamma(1 + n \beta)}, \ z \in C
$$
 is named as Mittag - Lefler function. \par

\vspace{3mm}

 Therefore, the function  $ s \to H_m(s), \ s \in S(m)  $   could be chosen as a density of distribution of a
random variable (vector), say, $ \kappa = \kappa_m  $ with support on the  simplex $  S(m): $

$$
{\bf P}(\kappa_m \in G) = \int_G H_m(s) \ ds \stackrel{def}{=} \mu_{m}(G).
$$

 This random vector $  \kappa = \kappa_{m} = \vec{\kappa}  = \vec{\kappa}_{m}   $ with values in
the polygon $  S(m) $   is a particular case of the so - called
{\it polygonal Beta distribution,} written:  $  \Law(\kappa) = PB(1/2,m), $   iff it has a
density $  H_{m}(s), \ s \in S(m),  $  see  e.g. \cite{Ostrovsky208}. \par

 Note that $ \dim \kappa_m = m.  $  The the most economical way
 of generation of this distribution, such that for each value
$ \kappa_m  $  is elapsed exact $ m $ random variables uniform distributed on the interval $ (0,1), $ is described in the
aforementioned article \cite{Ostrovsky208}. \par

\vspace{3mm}

We can offer the following probabilistic representations for the integral $ J_m(h): $

$$
  J_m(h) \cdot t^{-m/2}/W_m  = {\bf E} \tilde{h}(\zeta,\kappa_m),
$$
where the distribution  of the r.v. $  \kappa_m $ was described before, the random vector $ \zeta = \{ \zeta_1, \zeta_2, \ldots, \zeta_m   \} $
consists on the $  m $ independent centered Gaussian standard distributed {\it matrices} of the size $  d \times d: $

$$
f_{\zeta}(z_1,z_2,\ldots,z_m) = (2 \pi)^{-dm/2} \
\exp \left[ - \frac{||z_1||^2}{2}- \frac{||z_2||^2}{2} -  \ldots -  \frac{||z_m||^2}{2}  \right] \eqno(2.15)
$$
and the random vectors $ (\zeta, \kappa_m) $ are independent.\par

 In the sequel  the notation $  ||A||^2 $ for the $ d \times d $  matrix $ A = \{ a_{i,j}  \} $ with real entries denotes

$$
||A||^2 \stackrel{def}{=} \sum_{i=1}^d \sum_{j=1}^d a^2_{i,j}.
$$

  We denote also for the positive semi - definite matrix $  A  = \{ a_{i,j}  \}  $
  $$
  [A] := \max_i a_{i,i}.
  $$

 Note that for the r.v. $ (\zeta, \kappa_m) $  generation need $  m + dm = m(d+1) $ uniform distributed on the set $ [0,1] $ r.v.\par

\vspace{3mm}

 Let $  N(m) ( = N(m,n)  ) $ be arbitrary positive integer numbers. The Monte Carlo approximation $ J_{m, N(m)} $ for the integral
 $  J_m(h), $   the so-called "depending trial method", see \cite{Frolov1},  \cite{Grigorjeva1} has the form

$$
  J_{m, N(m)}(h) :=  t^{m/2} \cdot W_m  \cdot \frac{1}{N(m)} \sum_{i=1}^{N(m)}  \tilde{h}(\zeta_i,\kappa_{m,i}), \eqno(2.16)
$$
where the r.v. $  (\zeta_i, \kappa_{m,i} ) $ are independent copies of $ (\zeta, \kappa_m ). $ \par

Of course, this estimate is unbiased: $ {\bf E} J_{m, N(m)} = J_m(h). $ Let us estimate the variation of $  J_{m,N(m)}. $ Evidently,

$$
[\Var]{J_{m, N(m)}}(h) \le t^{m} \cdot W_m \cdot \frac{1}{N(m)} \cdot  J_m(\tilde{h}^2), \eqno(2.17)
$$
and following

$$
[\Var]{J_{m, N(m)}}(h) \le t^{m} \cdot W^2_m \cdot \frac{1}{N(m)} \cdot \sup_{y,s} h^2(y,s).  \eqno(2.17a)
$$

 Note that the integral $ J_m(\tilde{h}^2) $ may be estimated as well as the source integral $ J_m(h). $ \par

\vspace{4mm}

{\bf B. We denote and consider now the following integral }

$$
 I_m(h)  =  w \odot w \odot \ldots \odot  w \odot h =
$$

$$
\int_0^t ds_1 \int_0^{s_1} ds_2 \ldots \int_0^{s_{m-1}} ds_m \int_{R^{d m}} (2 \pi)^{-dm/2} (t-s_1)^{-d/2} (s_1 - s_2)^{-d/2} \ldots
(s_{m-1} - s_m)^{-d/2} \times
$$

$$
\exp \left[ - \frac{||x-y_1||^2}{2(t-s_1)}
- \frac{||y_1 - y_2||^2}{2(s_1 - s_2)} -  \ldots -  \frac{||y_{m-1} - y_m||^2}{2(s_{m-1} - s_m)}  \right] \times h(y,s) \ dy. \eqno(2.18)
$$

We have $  I_m(h) \cdot t^{-m} (2 \pi)^{dm/2} = $

$$
\int_{S(m)}  d \tau \int_{R^{d m}} \ dz  \times \exp \left[ - \frac{||z_1||^2}{2}
- \frac{||z_2||^2}{2} -  \ldots -  \frac{||z_m||^2}{2}  \right] \cdot \tilde{h}(z,\tau), \eqno(2.19)
$$
where  as before $  S(m)  $ is an  $ m \ - $  dimensional unit polygon  (simplex), $  (y,s) = Y(z,\tau) $ is at the
same substitution  and $  (z, \tau) = Y(y,s), \ \tilde{h}(z,\tau) = h(Y(y,s)) .  $\par

 Since the volume of the simplex $ S(m) $  is equal to $  1/m!, $ the expression for the integral $ I_m(h)  $ obeys a following
probabilistic representation:

$$
t^{-m} \cdot m! \cdot I_m(h) = {\bf E} \tilde{h}(\zeta, \nu_m), \eqno(2.20)
$$
where the r.v. $ \zeta  $ has as before multidimensional centered standard Gaussian (Normal) distribution in the space $ R^{d \ m}, $
the random vector $  \nu_m $ has an {\it  uniform  } distribution in the polygon $ S(m) $ and the r.v. $ (\zeta, \nu_m)  $ are independent. \par
  Recall that this distribution is a particular case of general Polygonal Beta distribution, see \cite{Ostrovsky208}. \par

 Naturally, the Monte Carlo approximation for $  I_m(h)  $ has a form

$$
I_{m, N(m)}(h):= t^m \cdot \frac{1}{m!} \cdot \frac{1}{N(m)} \ \sum_{i=1}^{N(m)}  \tilde{h}(\zeta_i, \nu_{m,i} ), \eqno(2.21)
$$
where $  N(m) $ is non-random positive integer  number, $ (\zeta_i, \nu_{m,i})  $ are independent copies of the r.v. $ (\zeta, \nu_m). $\par

This estimate is unbiased: $ {\bf E} I_{m, N(m)} = I_m(h). $ Let us estimate the variation of the random value $  I_{m,N(m)}. $

$$
[\Var]{I_{m, N(m)}}(h) \le t^{2m} \cdot \frac{1}{m!} \cdot \frac{1}{N(m)} \cdot I_m(\tilde{h}^2)  \eqno(2.22)
$$
and following

$$
[\Var]{I_{m, N(m)}}(h) \le t^{2m} \cdot \frac{1}{m!^2} \cdot \frac{1}{N(m)} \cdot \sup_{y,s} h^2(y,s).  \eqno(2.22a)
$$

 Note that the integral $ I_m(\tilde{h}^2) $ may be estimated as before as well as the source integral $ I_m(h). $ \par

\vspace{4mm}

{\bf  C. Mixed case.}\\

\vspace{3mm}

 We introduce for simplicity two operators:

$$
T_w[h](x,t) := [w \odot g](x,t),  \hspace{6mm}  T_{\partial w}[g](x,t) := [\partial w \odot g](x,t), \eqno(2.23)
$$
and consider the following multiple convolution:

$$
K[h] = K[h](y,s) = K_{m_1, m_2}[h](y,s) := T_w^{m_1} T_{\partial w}^{m_2}[h](y,s). \eqno(2.24)
$$

 Here $ m_1, m_2 $ be non - negative integer numbers  and we denote (temporarily) $ m = m_1 + m_2,  $ and suppose $ m \ge 1; $
the case $ m = 0 $ is trivial. \par
 We offer and investigate the Monte - Carlo approximation for $ K[h] $ computation alike in the pilcrow {\bf A} and in the pilcrow {\bf B.}\par

Figuratively speaking, the pilcrow {\bf C} is synthesis of the subsections {\bf A} and {\bf B.}\par

 We obtain after at the same linear changing of variables the following expression for the function $ K(\cdot) $

$$
K_{m_1,m_2}[h] = (2 \pi)^{dm/2} \cdot t^{m_1/2 + m_2} \cdot \int_{S(m_1)} d \tau \int_{S(m_2)} d\theta \int_{R^{d m}} dz \ R_m(\tau) \times
$$

$$
\exp \left( -0.5 \sum_{i=1}^m ||z_i||^2    \right)  \cdot \tilde{h}(z,\tau,\theta),
$$
which admit the next probabilistic representation

$$
K_{m_1,m_2}[h] = t^{m_1/2 + m_2} \cdot W(m_1) \cdot (1/m_2!) \cdot  {\bf E} \tilde{h} \left(\zeta,\kappa_{m_1}, \nu_{m_2} \right),  \eqno(2.25)
$$
where as before the random vectors $ ( \zeta, \kappa_{m_1}, \nu_{m_2}  ) $ are (common) independent, the random vector $ \zeta $
is normal centered standard blocky  distributed in the space $ R^{d \ m}, $ the r.v. $  \kappa_{m_1} $ has the Polygonal Beta distribution
with index  $ (-0.5) $ in the simplex $ S(m_1)  $ and the r.v. $ \nu_{m_2} $ is uniformly distributed inside the {\it other} simplex $  S(m_2). $\par

 The Monte  Carlo approximation  $ K_{m_1,m_2, N(m)}[h]  $ for the function $ K_{m_1, m_2}[h] $ is clear:

$$
 K_{m_1, m_2, N(m)}[h] := t^{m_1/2 + m_2} \cdot W(m_1) \cdot (1/m_2!) \cdot  \frac{1}{N(m)} \cdot
 \sum_{i=1}^{N(m)}  \tilde{h} \left(\zeta_i,\kappa_{m_1,i}, \nu_{m_2,i} \right),  \eqno(2.26)
$$
where  $ (\zeta_i,\kappa_{m_1,i}, \nu_{m_2,i}) $   are independent copies of $ (\zeta,\kappa_{m_1}, \nu_{m_2}). $ \par
 The approximation  $  K_{m_1, m_2, N(m)}[h] $ has the following variation estimate

 $$
 [\Var]( K_{m_1, m_2, N(m)}[h]) \le  t^{m_1 + 2 m_2} \cdot W^2(m_1) \cdot (1/m_2!^2) \cdot \sup_{y,s} h^2(y,s)/N(m), \eqno(2.27)
 $$
which allows in turn a very simple but rough estimate

 $$
 [\Var]( K_{m_1, m_2, N(m)}[h]) \le  \max(t^{m}, t^{ 2 m}) \cdot W^2(m) \cdot \sup_{y,s} h^2(y,s)/N(m). \eqno(2.28)
 $$

 Note that the number $  N(m) $ dependent only on the number  $ m, $ but not on the individual numbers $ m_1, \ m_2.$
Therefore, the general amount of standard uniformly distributed on the interval $ (0,1) $  random numbers for $ K_{m_1,m_2, N(m)}[h]  $
computations in (2.26) is equal to  $ N(m)\cdot (m + d m) = N(m)\cdot m \cdot (d+1). $  \par

\vspace{4mm}

\section{Auxiliary facts. Non-linear recursions.}

\vspace{4mm}

{\bf A. Polynomial recursion.} \\

\vspace{3mm}

 Let us introduce the following sequence of polynomials $ \{ P_n(z) \}, \ n = 0,1,2,\ldots: \ P_0(z) = 0, \  P_1(z)  = z $ and by recursion over $  n: $

$$
P_{n+1}(z) = z + P^2_n(z), \ n = 1,2,\ldots
$$

 For instance,

$$
P_2(z) = z + z^2 + 2 z^3 + z^4; \ P_3(z)= z + z^2 + 2z^3 + 5 z^4 + 6z^5 + 6z^6 + 4 z^7 + z^8.
$$

 Evidently, $  P_n(0) = 0 $ and   $ \deg(P_n) = 2^n. $ Further, we conclude by means of induction

$$
P_n(z) = \sum_{m=1}^{2^n} A(m,n) z^m,
$$
where the coefficients $  \{ A(m,n)  \} $ are positive and integer.\par
 We derive substituting the value $ z = 1 $ and denoting
 $$
  P_n(1) = \sum_{m=1}^{ 2^n} A(m,n) =: \tilde{P}(n)
 $$
the following recursion:
$$
\tilde{P}(0)= 0, \ \tilde{P}(1) = 1, \  \tilde{P}(n+1) = 1 + (\tilde{P}(n))^2.
$$

\vspace{3mm}

Suppose $ 0 < z < 1/4;  $ then the sequence $  \{ P_n(z) \}, \ n = 1,2,\ldots $ monotonically increases and converges uniformly
in the ball

$$
\{z: \ |z| \le 1/4 - \epsilon_0 \}, \ \epsilon_0 = \const \in (0,1/4)
$$
to the analytic function $  P(z) $  which satisfies the equation

$$
P(z) = 1 + z \ P^2(z), \ \Leftrightarrow P(z) = \frac{1- \sqrt{1-4z}}{2z}; \ P(0) = 1.
$$

 This function has a Taylor's expression

$$
P(z) = 1 +  \sum_{M=1}^{\infty} \frac{(2M)!}{ M! \ (M+1)!}  \ z^M.
$$

  As long as

$$
\frac{(2M)!}{ M! \ (M+1)!} \le \pi^{-1/2} \ M^{-1/2} \ (M+1)^{-1}  \ 2^{2M} \ e^{5/(24M)}, \ M \ge 1,
$$
the series for the function $ P(z) $  converge inside the closed ball $ |z| \le 1/4. $\par

 We  conclude also

 $$
\sup_n A(M,n) = \lim _{n \to \infty} A(M,n) = \frac{(2M)!}{ M! \ (M+1)!}.
 $$

 Moreover, let us introduce the following relation of equivalence  $ E(n), n = 2,3,\ldots, $  more exactly, the sequence of relations,
 between  the polynomials  $ P = P(z) $ and $ Q = Q(z): $

 $$
 P \stackrel{E(n)}{\sim } Q \Leftrightarrow \forall k \le n-1 \Rightarrow P = Q \mod z^k.
 $$
  We deduce then by induction:

 $$
  P_{n+1} \stackrel{E(n-1)}{\sim} P_n,
 $$
  and as a consequence

$$
\forall k \le n - 1 \Rightarrow A(k,n-1) = A(k,n).
$$

\vspace{3mm}

{\bf B.  Non - linear numerical recursion. }\par

\vspace{3mm}

 The following numerical recurrence relation (dynamical system) with quadratic non-linearity

$$
D(n+1) = 1 + d^2 \cdot D^2(n)
$$
with initial condition  $ D(0) = 1  $ is investigated with some numerical examples in the important case $  d = 3 $
 in particular in the article \cite{Ostrovsky209}. \par
 For instance, it was obtained there the speed of increase of solution $  D(n)  $ as $ n \to \infty $ and bilateral
exact bounds. \par
 For example:

 $$
 D(0) = 1, \ D(1) = 10, \ D(2) = 901, \ D(3) = 811 \ 802, \ D(4) = 659 \ 022 \ 487 \ 205,
 $$

$$
D(5) = 434 \ 310 \ 638 \ 641 \ 864 \ 388 \ 712 \ 026;
$$

\vspace{4mm}

$$
\forall k,l = 1,2,\ldots \ \Rightarrow
1 \le \frac{9 \ D(k + l)}{ [ 9 \ D(l)]^{2^k}  } \le \left[ 1 + \frac{1}{9 D^2(l)} \right]^{2^k - 1},
$$

$$
\forall k   \ge 1 \Rightarrow \lim_{l \to \infty} \frac{9\ D(k + l)}{ [ 9 \ D(l)]^{2^k}  } = 1.
$$

 It is proved also in \cite{Ostrovsky209} that $  D(n) $  is number of independent summands for the $ n^{th} $
iteration $ u_n, \ \partial u_n. $ \par

\vspace{3mm}

\section{Iterations.}

\vspace{3mm}

{\it Notations.} $ v = v_0 = u_0\cdot \partial u_0; $  \hspace{7mm}  $  |||u(\cdot, \cdot) ||| :=  $

$$
 \max_{i=1,2,\ldots,d} \max_{j=1,2,\ldots,d}   \sup_{x,t}
\max \{ | \max(|x|,1) \ u_i(x,t) \ |, \max(|x|,1) |\partial u_i/\partial x_j| \}. \eqno(4.0)
$$

\vspace{3mm}

 The mild solution $ u = u(x,t) $ of a Navier-Stokes equation in the whole space $  x \in R^d $ throughout
its lifetime $ t \in [0,T], T = \const \le \infty $ may be represented as a limit as $ n \to \infty, n = 0,1,2, \ldots $ the
following recursion:

 $$
 u_{n+1}(x,t) = u_0(x,t) + G[u_n, u_n](x,t), n=0,1,2,\ldots,
 $$
where $ u_0(x,t) $ is the solution of heat equation  with correspondent initial $  a(x) $ value and right-hand side $ f(x,t): $

$$
\partial u_0/\partial t = 0.5 \ \Delta u_0 + f(x,t), \ u_0(x, 0+) = a(x)
$$
and $ G[u,v] $ is bilinear unbounded pseudo-differential operator, \cite{Kato1}. \par

 The iterative recursion may be written as follows:

$$
u_{n+1} = u_0 + w \odot v_n, \hspace{6mm} \partial u_{n+1} = \partial  u_0 + \partial w \odot v_n,
$$
where

$$
v_n := u_n \cdot \partial u_n, \ n = 0,1,2,\ldots.
$$

 For instance,

$$
u_{1} = u_0 + w \odot v, \hspace{6mm} \partial u_{1} = \partial  u_0 + \partial w \odot v, \eqno(4.1)
$$

$$
u_{2} = u_0 + w \odot v_1, \hspace{6mm} \partial u_{2} = \partial  u_0 + \partial w \odot v_1,
$$
and we obtain after substitution

$$
u_2 = u_0 + w \odot v  + w \odot [ u_0 \cdot (\partial w \odot v ) ] +
$$

$$
w \odot \left\{ [w \odot v] \cdot \partial  u_0  \right\} + w \odot \left[(w \odot v ) \cdot ( \partial w \odot v) \right],\eqno(4.2a)
$$

$$
\partial u_2 = \partial u_0 + \partial w \odot v  + \partial w \odot [ u_0 \cdot (\partial w \odot v ) ] +
$$

$$
\partial w \odot \left\{ [w \odot v] \cdot \partial  u_0  \right\} + \partial w \odot \left[(w \odot v ) \cdot ( \partial w \odot v) \right].
\eqno(4.2b)
$$

\vspace{3mm}

 It follows from the inductions method that

 $$
 u_n = u_0 + \sum_{r=1}^{D(n)} Q_r(u_0), \eqno(4.3a)
 $$

 $$
 \partial u_n = \partial u_0 + \sum_{r=1}^{D(n)} \partial Q_r(u_0), \eqno(4.3b)
 $$
 where

$$
 Q_r(u_0) = Q_r(u_0, \partial u_0, v, \ldots)
$$
is {\it  homogeneous } relative the source function $ u_0(\cdot, \cdot) $  non - linear operator
acting on the  continuous differentiable functions defined on the semi - space $  R^d \times R_+  $ into itself.  \par

 Of course, we offer to compute each integral in (4.3a) and in (4.3b)  by means of the Monte Carlo method. We discuss many technical details.\par

 Note first of all that every member $ Q_r(\cdot) $ in (4.3a)  (and analogously in (4.3b) ) has a form

$$
Q_r(\cdot) = Q_{r; l_1,l_2}(\cdot) = Q^{(n)}_{r; l_1,l_2}(\cdot) = T_w^{l_1} T_{\partial w}^{l_2} [h],
$$
with appropriate function $ h = h[u_0] = h_{l_1, l_2}[u_0](y,s),  $  where

$$
l:= l_1 + l_2 = \deg_{u_0}(h)-1.\eqno(4.4)
$$

 The last notion   $ \deg_{u_0}[h] \stackrel{def}{=} k $ implies by definition that

$$
h[\lambda u_0] = \lambda^{k} h[u_0], \ \lambda = \const \in R.
$$

 The expression (4.3a) and equally (4.3b) may be rewritten as follows.

$$
u_n = u_0 + \sum_{k=1}^{2^n} L_k[u_0],
$$
where

$$
L_k[u_0] = L_k^{(n)}[u_0] = \sum \sum_{l_1,l_2: l_1 + l_2 = k} Q^{(n)}_{r; l_1,l_2}[h].\eqno(4.5)
$$

 Note that the amount of summands in the right - hand side of the expression (4.5) is equal to $  A(k,n). $\par
Each member in (4.5) has the degree  $  k-1 $ relative the function $ u_0(\cdot, \cdot) $ and may be computed
by means of the Monte Carlo method in accordance to the second section.\par

 We offer to give for all the members into (4.5) computation at the same amount $  N(k) = N(k,n)  $ random
(quasi - random) independent vector variables,  so that the {\it general} amount the standard (uniformly distributed) r.v.  (spending)
for approximation $  u_n \approx u_{n,N} $  Monte-Carlo computation is equal to

$$
N = \sum_{k=1}^{2^n} A(k,n) \cdot N(k,n)  \cdot d(k+1). \eqno(4.6)
$$

 Notice that at the same random variables may be  used also for the $ \partial u_n $ Monte  Carlo computation, for the sake of saving. \par

  We give now the rough variation estimate for $ u_{n,N}  $ approximation based on the formulae 2.28. Namely,

 $$
[ \Var]( u_{n,N}) \le \sum_{k=1}^{2^n}   \max(t^{k}, t^{ 2 k}) \cdot W^2(k) \cdot A(k,n) \cdot |||u_0|||^{ 2k } /N(k,n) \eqno(4.7)
 $$
and at the same estimate is true  for $ \partial u_n \approx \partial u_{n,N}  $ Monte Carlo computation. \par

\vspace{3mm}

{\bf Remark 4.1.} Recall that the functions $ u, \ \partial u  $ and following $ u_n, \ \partial u_n  $ are vector and moreover
matrix functions. For instance, $ \partial u_i^{(j)} = \partial u_i/\partial x_j.  $ \par

 But for the $ \partial u_i /\partial x_j $ by means of offered here method can be used, for the sake of saving,
at the same random variables  as by computation $ u_n = u_{n,i}. $ \par

\vspace{4mm}

\section{Subject of optimization}

\vspace{4mm}

It seems quite reasonable the following statement of constrained optimization problem. Let  the general amount of spending
standard distributed r.v. $  N  $ be a given "great" number, for example, $ N = 10^6 \ - \ 10^8.  $ \par

 Subject of minimization:

 $$
Z:= Z(N(1), N(2), \ldots, N(2^n)) \stackrel{def}{=} \sum_{k=1}^{2^n} W^2(k) \cdot A(k,n) \cdot |||u_0|||^{ 2k } /N(k) \eqno(5.1)
 $$

 This function is weakly proportional to the upper estimation for the variance  $ \Var( u_{n,N}) $ in (4.7), moreover:

 $$
[\Var]( u_{n,N})  \asymp  Z(N(1), N(2), \ldots, N(2^n)), \ t \in (0,T), \ T = \const \in (0,\infty).
 $$

Restriction:

$$
Y := Y(N(1), N(2), \ldots, N(2^n))\stackrel{def}{=} \sum_{k=1}^{2^n} A(k,n) \cdot N(k)  \cdot d(k+1) = N. \eqno(5.2)
$$

 So, we get to the following problem of constrained optimization:

$$
Z(N(1), N(2), \ldots, N(2^n)) \to \min  / Y(N(1), N(2), \ldots, N(2^n)) = N, \ N(k) \ge 1. \eqno(5.3)
$$

 We find by means of Lagrange factor method neglecting to start an integer variables:

$$
N_0(k) = \frac{N \ d^{-3/2} \ W(k) \ |||u_0|||^k}{\sum_{r=1}^{2^n} W(r) \ |||u_0|||^r},
$$
up to around to nearest integer number.   Herewith

$$
\min_{ \{ N(k)  \}} [\Var](u_{n,N}) \asymp \frac{d^{3/2}}{N} \cdot \sum_{m=1}^{2^n} W(m) \cdot A(m,n) \cdot |||u_0|||^m \times
$$

$$
\sum_{r=1}^{2^n} (r+1)^{1/2} \cdot W(r) \cdot A(r,n) \cdot ||| u_0|||^r. \eqno(5.4)
$$

 {\it Some slight simplification:} as $ N \to \infty $

$$
 \min_{ \{ N(k)  \}} [\Var](u_{n,N}) \asymp \frac{d^{3/2} C(n)}{N} \eqno(5.5a)
$$
and analogously

$$
 \min_{ \{ N(k)  \}}  [\Var]( \partial u_{n,N}) \asymp \frac{d^{3/2} C(n)}{N}.\eqno(5.5b)
$$

 The last estimates imply that the speed of convergence $  u_{n,N}   $ to $  u_n $ as $  N \to \infty $
is equal $  N^{-1/2}, $  as in the classical Monte Carlo method. \par
 For the linear integral equations this effect was detected in \cite{Frolov1}, \cite{Grigorjeva1}.\par

\vspace{4mm}

\section{Concluding remarks.}

\vspace{4mm}

{\bf A. Functional approach. } \par

\vspace{4mm}

 In order to estimate the  random error, i.e. the deviation $  u_{n,N} - u_n  $ in some space - time norm $  || \cdot|| = ||\cdot||X,T,  $
more exactly, to estimate the value

$$
Q_{n,N}( v ) \stackrel{def}{=} {\bf P} ( \sqrt{N} || u_{n,N} - u_n ||X,T  > v), \ n = 2,3,\ldots, \  v = \const > 0,  \eqno(6.0)
$$
we to use the Central Limit Theorem  (CLT) in the correspondent Banach space $ (K,||\cdot||X,T $ , in accordance with which there exists a limit

$$
\lim_{N \to \infty} Q_{n,N}( v )  \stackrel{def}{=} Q_{n,\infty}( v ),\eqno(6.1)
$$
where

$$
Q_{n,\infty}( v ) = {\bf P} (|| \xi(\cdot, \cdot) || > v),\eqno(6.2)
$$

$  \xi(x,t) $ is centered Gaussian distributed random field with values in the space $  K, \ (\mod {\bf P}). $ \par
 For the space of continuous functions  it is proved, e.g. in \cite{Frolov1}, \cite{Dudley1}, \cite{Gine1}, \cite{Grigorjeva1}, \cite{Heinkel1},
\cite{Kozachenko1}, \cite{Ostrovsky1}, \cite{Vaart1}; in the classical Lebesgue - Riesz  spaces $  L_p \ - $  in \cite{Ledoux1};
in the mixed Lebesgue-Riesz spaces-in \cite{Ostrovsky210},  \cite{Ostrovsky211}; in the mixed hybrid Lebesgue-continuous spaces-
in  \cite{Ostrovsky212} etc. \par

 The  behavior as $ v \to \infty $ of the probability $ Q_{n,\infty}( v ),  $ asymptotical as well as non - asymptotical is obtained in many works,
see, e.g. \cite{Piterbarg1}:

$$
Q_{n,\infty}( v ) \sim C(X,T) \ v^{\kappa-1} \ \exp(-v^2/(2 \sigma^2)), \ C(X,T), \kappa,  \sigma^2 = \const > 0. \eqno(6.3)
$$

 Equating the approximation of a value $ Q_{n,\infty}( v ) $  with $ v \ge 3 \sigma, $
 in the right - hand  of  (6.3) to the value $  \delta: $

$$
C(X,T) \ v(\delta)^{\kappa-1} \ \exp(-v(\delta)^2/(2 \sigma^2)) = \delta, \ v(\delta) \ge 3 \sigma, \eqno(6.4)
$$
 where $  1- \delta $ is reliability of the confidence interval, for instance, $ 0.95 $ or $ 0.99, $ we obtain the  asymptotical
confidence region for  the function $ u_n(x,t) $ in the norm $ || \cdot|| $ of the form

$$
|| u_n - u_{n,N}|| \le \frac{v(\delta)}{\sqrt{N}}.\eqno(6.5)
$$

\vspace{3mm}

{\bf B. General optimization.}\par

\vspace{3mm}

 The inequality (1.0a)  follows from  (6.2) - (6.5). Moreover,  it may be proved under simple condition that
there exists finite  function $  K = K_p(n), \ p > 1 $ for which

$$
\forall n = 1,2,\ldots  \ \Rightarrow  \left[ {\bf E} ||  u_{n,N}  - u_n||^p \right]^{1/p} \le K_p(n) \ N^{-1/2},
$$
and analogous conclusion may be obtained for the $ G(\psi) $  norm for the norm  difference  $  ||  u_{n,N}  - u_n||.$ \par

 The accuracy calculation show us that the constant $  K(n) $
in (1.0a) is proportional to the value $  D(n),  $ where $  D(n) $ is introduced and investigated in the third section, and obviously the relation
$  \sqrt{N} >> D(n) $ should be performed.\par
 The common error $ ||u - u_{n,N}||, $  including the deterministic part $ \le q^n $ and probabilistic part $ \le D(n)/\sqrt{N} $  does not exceed
the value

$$
\Delta := C \left( q^n + \frac{D(n)}{\sqrt{N}} \right). \eqno(6.6)
$$
 It appears naturally the   following qualitative optimization problem by fixed great value $ N:  $

$$
q^n + \frac{D(n)}{\sqrt{N}} \to \min_n: \ D(n)  << \sqrt{N}. \eqno(6.7)
$$
 The practical computation  taking into account the rate of increasing of the sequence $  \{  D(n) \} $
 show us that the optimal value $ n $ is 4 - 5. \par

\end{document}